\documentclass[10pt]{article}
\usepackage{mathrsfs}
\usepackage[all]{xy}
\usepackage{amssymb,latexsym,amsmath,amsthm, mathrsfs}
\usepackage{graphicx}
\usepackage{setspace}
\usepackage{fullpage}

\newtheorem{theorem}{Theorem}[section]
\newtheorem{lemma}[theorem]{Lemma}

\newcommand{\susp}{\Sigma\mkern-10.1mu/}

\newtheorem{definition}[theorem]{Definition}
\newtheorem{example}[theorem]{Example}

\newtheorem{proposition}[theorem]{Proposition}
\newtheorem{corollary}[theorem]{Corollary}
\newenvironment{acknowledgements}[1][Acknowledgements]
{\begin{trivlist} \item[\hskip \labelsep {\bfseries #1}]}
{\end{trivlist}}
\newenvironment{proof*}[1][proof]
{\begin{trivlist} \item[\hskip \labelsep {\bfseries #1}]}
{\end{trivlist}}

\newtheorem{conjecture}[theorem]{Conjecture}
\newtheorem{remark}[theorem]{Remark}
 \DeclareMathSymbol{\N}{\mathbin}{AMSb}{"4E}

\DeclareMathSymbol{\Z}{\mathbin}{AMSb}{"5A}
\DeclareMathSymbol{\R}{\mathbin}{AMSb}{"52}
\DeclareMathSymbol{\Q}{\mathbin}{AMSb}{"51}
\DeclareMathSymbol{\I}{\mathbin}{AMSb}{"49}
\DeclareMathSymbol{\C}{\mathbin}{AMSb}{"43}

\def\P{{\mathbb P}}

\def\A{{\mathbb A}}

\def\L{{\mathbb L}}

\numberwithin{equation}{section}

\title{Motivic integration and projective bundle theorem in morphic cohomology}

\author{Jyh-Haur Teh}

\begin{document}
\maketitle

\begin{abstract}
We reformulate the construction of Kontsevich's completion and use
Lawson homology to define many new motivic invariants. We show that
the dimensions of subspaces generated by algebraic cycles of the
cohomology groups of two $K$-equivalent varieties are the same,
which implies that several conjectures of algebraic cycles are
$K$-statements. We define stringy functions which enable us to ask
stringy Grothendieck standard conjecture and stringy Hodge
conjecture. We prove a projective bundle theorem in morphic
cohomology for trivial bundles over any normal quasi-projective
varieties.
\end{abstract}

\section{Introduction}
With the insight from string theory, Batyrev first showed that two
birational Calabi-Yau manifolds have the same Betti numbers by using
Weil's conjecture (see \cite{Baty2}). Kontsevich pushed this result
a step further by introducing the notion of motivic integration in
showing that two $K$-equivalent varieties have the same Hodge
numbers. In this paper we show further that the dimensions of
subspaces generated by algebraic cycles of cohomology groups of two
$K$-equivalent varieties are the same. We found that to construct
motivic measure and motivic integration, we do not need the product
structure of $K_0(Var)$, the Grothendieck group of algebraic
varieties. It is sufficient to give $K_0(Var)$ a
$\mathscr{L}$-module structure where $\mathscr{L}$ is the polynomial
ring generated by $\L=\C$, and then we can follow Kontsevich's
construction to get an abelian group completion. This makes a huge
difference since many invariants are not multiplicative, but they
can be defined in our new settings. The tools we need to use are the
homotopy property and the blow-up formula in Lawson homology. We
first show that we can define invariants over the Grothendieck group
$K_0(SPV)$ of smooth projective varieties, and then show that the
isomorphism given by Bittner $\varphi: K_0(Var)\rightarrow K_0(SPV)$
is an isomorphism of $\mathscr{L}$-modules. Through this
isomorphism, invariants defined on smooth projective manifolds
induce invariants for quasi-projective varieties. These invariants
induce invariants defined on the image of some localization
$S^{-1}K_0(Var)$ of $K_0(Var)$ in the Kontsevich's completion which
enables us to use motivic integration.

We review and modify some constructions of motivic measure and
motivic integration in Section 2, use the blow-up formula in Lawson
homology and the natural transformations from Lawson homology to
singular homology to define some motivic invariants in Section 3.
This enables us to show that the generalized Hodge conjecture, the
Grothendieck standard conjecture, the Friedlander-Lawson conjecture,
the Friedlander-Mazur conjecture are $K$-statements. In Section 4 we
define stringy functions which extend many classical notions to
varieties with singularities. One of our most interest stringy
functions defines the stringy version of the dimension of cohomology
classes generated by algebraic cycles for singular varieties. We are
then able to ask stringy Grothendieck standard conjecture and
stringy Hodge conjecture. The stringy Grothendieck standard
conjecture is verified for normal projective toric varieties with
$\Q$-Gorenstein singularities. We conjecture that for mirror pairs
$(V, W)$ of dimension $n$ constructed by Batyrev and Borisov, the
relation of Hodge numbers $h^{p, q}(V)=h^{n-p, q}(W)$ can be
enhanced to a relation of stringy $\phi$-numbers.

In Section 5 we focus on varieties with finitely generated Lawson
homology groups. We show that motivic integration can be defined
over these varieties, and we do the same thing for higher Chow
groups. Since one of the main tools we use in this paper is the
projective bundle theorem in Lawson homology, it is natural to ask
if similar result holds in morphic cohomology. Friedlander proved a
projective bundle theorem (\cite{F4}) in morphic cohomology for
smooth quasi-projective varieties but since there is no
Mayer-Vietoris sequence in morphic cohomology at this moment, a
proof of the result for general quasi-projective varieties is
difficult to get. In section 6, we are able to prove a projective
bundle theorem in morphic cohomology for trivial bundles over any
normal quasi-projective varieties without assuming smoothness. This
seeming trivial result already applies almost all techniques in
Lawson homology and morphic cohomology.

\section{Motivic integration}
\subsection{Arc spaces}
Let $SPV$ be the collection of all isomorphism classes of smooth
projective varieties and $Var$ be the collection of all
quasi-projective varieties. Let $K_0(SPV)=\Z(SPV)/\sim_{bl}$ be the
free abelian group generated by elements in $SPV$ quotient by the
subgroup $\sim_{bl}$ which is generated by elements of the form
$Bl_YX-X+Y-E(Y)$ where $Y$ is a smooth subvariety of a smooth
projective variety $X$, $Bl_YX$ is the blow-up of $X$ along $Y$ and
$E(Y)$ is the exceptional divisor of this blow-up. Let
$K_0(Var)=\Z(Var)/\sim$ be the Grothendieck group of
quasi-projective varieties. The subgroup $\sim$ is generated by
elements of the form $X-(X\backslash Y)-Y$ where $Y$ is a locally
closed subvariety of $X$.

Let $\L=\C$ and let $\mathscr{L}=\Z\{\L^i|i\in \Z_{\geq 0}\}$ be the
free abelian group generated by $\L^i$ for all nonnegative integer
$i$. $\mathscr{L}$ is a ring with the obvious multiplication. For
$[X]\in K_0(SPV)$, we define $\L^0\cdot [X]=[X]$ and
$$\L^i\cdot [X]:=[X\times \P^i]-[X\times \P^{i-1}]$$ for $i>0$.
Since
$$[(Bl_YX-X+Y-E(Y))\times \P^i]=[Bl_{Y\times \P^i}(X\times
\P^i)-X\times \P^i+Y\times \P^i-E(Y\times \P^i)]$$ this
multiplication is well defined on $K_0(SPV)$, and it is easy to
check that $\L^i(\L^j[X])=\L^{i+j}[X]$, the group $K_0(SPV)$ becomes
a $\mathscr{L}$-module under this action. The group $K_0(Var)$ is
naturally a $\mathscr{L}$-module under the product of varieties.

Let $S=\{\L^i\}_{i\geq 0}\subset \mathscr{L}$ and let
$\mathcal{M}=S^{-1}K_0(SPV)=\{\frac{a}{\L^i}|a\in K_0(SPV), i\in
\Z_{\geq 0}\}$, $\mathcal{N}=S^{-1}K_0(Var)=\{\frac{a}{\L^i}|a\in
K_0(Var), i\in \Z_{\geq 0}\}$ be the group obtained by taking the
localization of $K_0(SPV), K_0(Var)$ respectively with respect to
the multiplicative set $S$. Define $F^k\mathcal{M}, F^k\mathcal{N}$
to be the subgroups of $\mathcal{M}, \mathcal{N}$ generated by
elements of the form $\frac{[X]}{\L^i}$ where $i-dim X\geq k$. Then
we get a decreasing filtration
$$\cdots \supseteq F^k\mathcal{M} \supseteq F^{k+1}\mathcal{M} \supseteq \cdots $$
of abelian subgroups of $\mathcal{M}$ and a decreasing filtration
$$\cdots \supseteq F^k\mathcal{N} \supseteq F^{k+1}\mathcal{N} \supseteq \cdots $$
of abelian subgroups of $\mathcal{N}$.

Suppose that we are given a decreasing filtration $ \cdots \supseteq
F^k \supseteq F^{k+1} \supseteq \cdots $. A Cauchy sequence with
respect to this filtration is a sequence $\{a_i\}$ where $a_i\in
F^i$ for all $i$ such that for any $n>0$, there is a $N>0$ such
$a_i-a_j\in F^n$ for all $i, j>N$.

\begin{definition}
The Kontsevich group of smooth projective varieties is defined to be
$$\hat{\mathcal{M}}:=\underset{\leftarrow}{\lim}\frac{\mathcal{M}}{F^k\mathcal{M}}$$
the completion of $\mathcal{M}$ with respect to the filtration
above. Similarly, we define the Kontsevich group of varieties to be
$$\hat{\mathcal{N}}:=\underset{\leftarrow}{\lim}\frac{\mathcal{N}}{F^k\mathcal{N}}$$
the completion of $\mathcal{N}$ with respect to the filtration of
$\mathcal{N}$.

We use also $F^{\bullet}\mathcal{M}, F^{\bullet}\mathcal{N}$ to
denote the filtrations in $\hat{\mathcal{M}}, \hat{\mathcal{N}}$
respectively induced by the filtrations above.  We denote
$\overline{\mathcal{M}}$ to be the image under the canonical map
$\mathcal{M}\rightarrow \hat{\mathcal{M}}$, and
$\overline{\mathcal{N}}$ to be the image under the canonical map
$\mathcal{N}\rightarrow \hat{\mathcal{N}}$.
\end{definition}

\begin{definition}
We give $\hat{\mathcal{N}}$ a $\mathscr{L}$-module structure as
following: for a Cauchy sequence $(a_1, a_2, ...)\in
\hat{\mathcal{N}}$, define
$$\L\cdot (a_1, a_2, ....):=(\L a_1, \L a_2, ...)$$ which is again a
Cauchy sequence. It is easy to see that the canonical map
$\phi:\mathcal{N} \rightarrow \hat{\mathcal{N}}$ defined by
$\phi(a)=(a, a, ...)$ is a morphism of $\mathscr{L}$-modules and
$\overline{\mathcal{N}}$ is a submodule of $\hat{\mathcal{N}}$. We
define the $\mathscr{L}$-module structure similarly for
$\hat{\mathcal{M}}$.
\end{definition}

Let us recall a result of Bittner \cite{Bit}. For a better
presentation of the proof see \cite{MP}.
\begin{theorem}
There is a group isomorphism $\varphi: K_0(Var) \rightarrow
K_0(SPV)$.
\end{theorem}

The isomorphism $\varphi$ is given inductively on the dimension of
varieties. Assume it is defined for varieties of dimension less than
$n$. If $dim X=n$, we consider two cases:
\begin{enumerate}
\item If $X$ is nonsingular, let $\overline{X}$ be a nonsingular
compactification of $X$, then define
$$\varphi(X):=\overline{X}-\varphi(\overline{X}-X)$$

\item If $X$ is singular, let $X=\coprod_i S_i$ be a stratification of
$X$, then define
$$\varphi(X):=\sum_i \varphi(S_i)$$
\end{enumerate}

\begin{proposition}
\begin{enumerate}
\item The isomorphism $\varphi: K_0(Var) \rightarrow
K_0(SPV)$ induces an isomorphism $\varphi: \mathcal{N} \rightarrow
\mathcal{M}$ of $\mathscr{L}$-modules.
\item It induces an isomorphism of $\mathscr{L}$-modules $\hat{\varphi}:
\hat{\mathcal{N}} \rightarrow \hat{\mathcal{M}}$.
\item It induces an isomorphism of $\mathscr{L}$-modules $\overline{\varphi}:
\overline{\mathcal{N}}\rightarrow
\overline{\mathcal{M}}$.\label{isomorphism}
\end{enumerate}
\end{proposition}

\begin{proof*}
\begin{enumerate}
\item It suffices to prove $\varphi(\L \cdot X)=\L\cdot \varphi(X)$ for
$X$ nonsingular. It is easy to check for $X$ of dimension 1. We
assume that it is true for varieties of dimension less than $n$.
Then for $dim X=n$, we have $\varphi(\L\cdot X)=\overline{X\times
\C}-\varphi(\overline{X\times \C}-X\times \C)=\overline{X}\times
\P^1-\varphi((\overline{X}-X)\times
\C+\overline{X})=\overline{X}\times
\P^1-\overline{X}-\varphi((\overline{X}-X)\times \C)=\L\cdot
\overline{X}-\L\cdot \varphi(\overline{X}-X)=\L\cdot \varphi(X)$.

\item For $\frac{[X]}{\L^i}\in \mathcal{M}$, we define
$\hat{\varphi}(\frac{[X]}{\L^i})=\frac{\varphi(X)}{\L^i}$. Hence
$\hat{\varphi}(F^k\mathcal{N})\subset F^k\mathcal{M}$, and we have
$\hat{\varphi}^{-1}(F^k\mathcal{M})\subset F^k\mathcal{N}$.
Therefore $\varphi$ induces an isomorphism
$\hat{\varphi}:\hat{\mathcal{N}} \rightarrow \hat{\mathcal{M}}$ on
the completions. This is obviously an isomorphism of
$\mathscr{L}$-modules since $\varphi$ is an isomorphism of
$\mathscr{L}$-modules.

\item We have a commutative diagram
$$\xymatrix{ \mathcal{N} \ar[r]^{\varphi} \ar[d]& \mathcal{M} \ar[d]\\
\hat{\mathcal{N}} \ar[r]^{\hat{\varphi}} & \hat{\mathcal{M}}}$$
where $\varphi, \hat{\varphi}$ are $\mathscr{L}$-module isomorphisms
which implies that we have a $\mathscr{L}$-module isomorphism
between $\overline{\mathcal{N}}$ and $\overline{\mathcal{M}}$.
\end{enumerate}
\end{proof*}

From this Proposition, once we have a group homomorphism from
$\overline{\mathcal{M}}$ to some group $G$, we can use
$\overline{\varphi}$ to define a group homomorphism from
$\overline{\mathcal{N}}$ to $G$ which means that we can define an
invariant for all quasi-projective varieties.

\subsection{Motivic integration}
We give a brief review of the arc spaces of quasi-projective
varieties here. For the details, we refer to \cite{Baty} and
\cite{DL2}. We work over the field of complex numbers. For a complex
algebraic variety $X$ of dimension $d$, the space of $n$-arcs on $X$
is defined to be
$$\mathcal{L}_n(X)=Mor_{\C-\mbox{schemes}}(Spec\C[[t]]/(t^{n+1}), X).$$
For $m\geq n$, there are canonical morphisms $\theta^m_n:
\mathcal{L}_m(X) \rightarrow \mathcal{L}_n(X)$. Taking the
projective limit of these algebraic varieties $\mathcal{L}_n(X)$, we
obtain the arc space $\mathcal{L}(X)$ of $X$. For every $n$ we have
a natural morphism
$$\pi_n: \mathcal{L}(X) \rightarrow \mathcal{L}_n(X)$$
obtained by truncation. A subset $A$ of $\mathcal{L}(X)$ is called
cylindrical if $A=\pi^{-1}_n(C)$ for some $n$ and some constructible
subset $C$ of $\mathcal{L}_n(X)$. We say that $A$ is stable at level
$n$ if furthermore the restriction of
$\pi_{m+1}(\mathcal{L}(X))\rightarrow \pi_m(\mathcal{L}(X))$ over
$\pi_m(A)$ is a piecewise Zariski fibration over $\pi_m(A)$ with
fiber $\C^d$ for all $m\geq n$. We call $A$ stable if it is stable
at some level $n$. If $X$ is smooth, then all cylindrical sets of
$\mathcal{L}(X)$ are stable.

In the following, let us recall some constructions and results in
motivic integration. Even though we have almost all the
constructions and results from classical motivic integration, we
note that we only consider $K_0(Var), \overline{\mathcal{N}},
\hat{\mathcal{N}}$ as $\mathscr{L}$-modules, not rings.

\begin{definition}
If $A$ is stable at level $n$, we define
$$\widetilde{\mu}(A)=[\pi_n(A)]\L^{-nd}$$
in $\mathcal{N}$. Let
$$\mathcal{L}^{(e)}(X):=\mathcal{L}(X)\backslash
\pi^{-1}_e(\pi_e(\mathcal{L}(X_{sing})))$$ where $X_{sing}$ denote
the singular locus of $X$ and we view $\mathcal{L}(X_{sing})$ as a
subset of $\mathcal{L}(X)$. For a cylindrical set $A$, it can be
proved that $A\cap \mathcal{L}^{(e)}(X)$ is stable and we define
$$\mu(A)=\lim_{e \to \infty}\widetilde{\mu}(A \cap
\mathcal{L}^{(e)}(X))\in \hat{\mathcal{N}}$$
\end{definition}

Define a norm $||\cdot||$ on $\hat{\mathcal{N}}$ by $||a||:=2^{-n}$
where $n$ is the largest $n$ such that $a\in F^n\mathcal{N}$.

Then
\begin{enumerate}
\item for all $a, b \in \hat{\mathcal{N}}$, $||a+b||\leq max(||a||,
||b||)$,

\item for any $A, B$ cylindrical sets, we have $||\mu(A\cup B)||\leq
max(||\mu(A)||, ||\mu(B)||)$ and $||\mu(A)||\leq ||\mu(B)||$ when
$A\subset B$.
\end{enumerate}

\begin{definition}
We say that a subset $A$ of $\mathcal{L}(X)$ is measurable if, for
every positive real number $\epsilon$, there exists a sequence of
cylindrical subsets $A_i(\epsilon)$, $i\in \N$ such that
$$(A\Delta A_0(\epsilon))\subset \underset{i\geq 1}{\cup}A_i(\epsilon)$$
and $||\mu(A_i(\epsilon))||\leq \epsilon$ for all $i\geq 1$. We say
that $A$ is strongly measurable if moreover we can take
$A_0(\epsilon)\subset A$.
\end{definition}

The following is the result A.6 from \cite{DL3}.
\begin{theorem}
If $A$ is a measurable subset of $\mathcal{L}(X)$, then
$$\mu(A):=\underset{\epsilon \to 0}{\lim}\mu(A_0(\epsilon))$$
exists in $\hat{\mathcal{N}}$ and is independent of the choice of
the sequences $A_i(\epsilon), i\in \N$.
\end{theorem}

\begin{definition}
Let $X$ be a quasi-projective variety of pure dimension $d$. We
define the motivic volume of $X$ to be $\mu(\mathcal{L}(X))\in
\hat{\mathcal{N}}$. It can be shown that
$$\mu(\mathcal{L}(X))=\lim_{n\to
\infty}\frac{[\pi_n(\mathcal{L}(X))]}{\L^{nd}}$$ and it equals to
$[X]$ when $X$ is nonsingular.
\end{definition}

\begin{definition}
Let $A \subset \mathcal{L}(X)$ be measurable and $\alpha:A
\rightarrow \Z\cup \{\infty\}$ a function with measurable fibres
$\alpha^{-1}(n)$ for $n\in \Z$. We define the motivic integral of
$\alpha$ to be
$$\int_A\L^{-\alpha}d\mu:=\sum_{n\in \Z}\mu(\alpha^{-1}(n))\L^{-n}$$
in $\hat{\mathcal{N}}$ whenever the right hand side converges in
$\hat{\mathcal{N}}$, in which case we say that $\L^{-\alpha}$ is
integrable on $A$. If $\alpha$ is bounded from below, this is always
the case (see \cite{DL2}).
\end{definition}

\begin{definition}
Let $\mathcal{I}$ be a sheaf of ideals on $X$. We define
$$ord_t\mathcal{I}:\mathcal{L}(X) \rightarrow \N\cup \{\infty\}$$
by $ord_t\mathcal{I}(\gamma)=\underset{g}{min}\{ord_tg(\gamma)\}$
where the minimum is taken over $g\in \mathcal{I}$ in a neighborhood
of $\pi_0(\gamma)$. For an effective Cartier divisor $D$, we define
$ord_tD=ord_tI$ where $I$ is the ideal sheaf associated to $D$.
\end{definition}

The following result is from Theorem 2.7.1 of \cite{Loe1}.
\begin{theorem}(Change of variables formula)
Let $X$ be a complex algebraic variety of dimension $d$. Let $h: Y
\rightarrow X$ be a proper birational morphism and $Y$ a smooth
variety. Let $A$ be a subset of $\mathcal{L}(X)$ such that $A$ and
$h^{-1}(A)$ are strongly measurable. Assume that $\L^{-\alpha}$ is
integrable on $A$. Then
$$\int_A\L^{-\alpha}d\mu=\int_{h^{-1}(A)}\L^{-\alpha\circ h-ord_th^*(\Omega^d_X)}d\mu$$
where $h^*(\Omega^d_X)$ is the pullback of the sheaf of regular
differential $d$-forms of $X$.
\end{theorem}

For a divisor $D=\sum^r_{i=1}a_iD_i$ on $X$ and any subset $J\subset
\{1,..., r\}$, denote
$$D_J=\left\{
    \begin{array}{ll}
      \underset{j\in J}{\cap}D_j, & \hbox{ if } J\neq \emptyset \\
      X, & \hbox{ if } J=\emptyset.
    \end{array}
  \right.
$$
and $D^0_J:=D_J-\cup_{i\notin J}D_i$.

Even though we do not have $[X\times Y]=[X][Y]$ in
$\hat{\mathcal{N}}$ for general varieties $X, Y$, we do have
$[Y\times \C^*]=[Y\times (\C-0)]=[Y\times \C]-[Y\times 0]=\L
[Y]-[Y]=[Y][\C^*]$. By a similar calculation as in Theorem 6.28 of
\cite{Baty}, we have the following result.
\begin{theorem}
Let $X$ be a nonsingular algebraic varieties of dimension $d$ and
$D=\sum^r_{i=1}a_iD_i$ an effective divisor on $X$ with only simple
normal crossings. Then
$$\int_{\mathcal{L}(X)}\L^{-ord_tD}d\mu=\underset{J\subset \{1,...,
r\}}{\sum}[D^0_J]\cdot(\underset{j\in
J}{\prod}\frac{\L-1}{\L^{a_j+1}-1})$$ where $J$ is any subset
(including empty set) of $\{1,..., r\}$. \label{formula}
\end{theorem}

\begin{corollary}
Let $X$ be a variety of pure dimension $d$, and let $h:Y \rightarrow
X$ be a resolution of singularities of $X$ such that the relative
canonical divisor $D=\sum^r_{i=1}a_iD_i=K_Y-h^*K_X$ of $h$ has
simple normal crossings. Furthermore, assume that the ideal sheaf
$h^*(\Omega^d_X)$ is invertible. Then
$$\mu(\mathcal{L}(X))=\underset{J\subset \{1,..., r\}}{\sum}[D^0_J]\cdot(\underset{j\in
J}{\prod}\frac{\L-1}{\L^{a_j+1}-1}).$$ Hence $\mu(\mathcal{L}(X))$
belongs to $\overline{\mathcal{M}}[(\frac{1}{\L^i-1})_{i \geq 1}]$
\end{corollary}

We say that two smooth projective varieties $X$ and $Y$ are
$K$-equivalent if there is a smooth projective variety $Z$ and
birational morphisms $\rho_1:Z \rightarrow X$ and $\rho_2:Z
\rightarrow Y$ such that $\rho_1^*K_X=\rho_2^*K_Y$ in $Z$ where
$K_X, K_Y$ are the canonical divisors on $X$ and $Y$ respectively.
As a simple consequence of the ``change of variables formula", we
have the following result.

\begin{theorem}
If two smooth projective varieties $X, Y$ are $K$-equivalent, then
$[X]=[Y]$ in $\overline{\mathcal{N}}[(\frac{1}{\L^i-1})_{i \geq
1}]$. \label{$K$-equivalent}
\end{theorem}

\section{Lawson homology groups}
For an overview of Lawson homology and morphic cohomology, we refer
to \cite{FL3, FL1}. Recall that in Lawson homology we have the
homotopy property $L_pH_n(X\times \C^t)=L_{p-t}H_{n-2t}(X)$ for $X$
a smooth projective variety.

\begin{definition}
For a quasi-projective variety $X$ and an integer $p<0$, we define
the negative cycle group $Z_p(X):=Z_0(X\times \C^{-p})$, and
$H_p(X):=H^{BM}_0(X\times \C^{-p})$. Then we have
$L_pH_k(X)=L_0H_{k-2p}(X\times \C^{-p})=H^{BM}_{k-2p}(X\times
\C^{-p})$. Throughout this paper we will identify $L_rH_k(X)$ with
$L_{r+t}H_{k+2t}(X\times \C^t)$.
\end{definition}

Recall that we have natural transformations $\Phi:L_pH_k(X)
\rightarrow H_k(X)$ from Lawson homology to singular homology (see
\cite{F1}, \cite{FL1}, \cite{L1}). The intersection theory of cycle
spaces was developed by Friedlander and Gabber in \cite{FG} in which
they obtained a projective bundle theorem for Lawson homology. We
extend their result to negative cycle groups.

\begin{proposition}
Let $E$ be an algebraic vector bundle of rank $r+1$ over a
quasi-projective variety $Y$ of dimension $n$ and $\pi:
P(E)\rightarrow Y$ the projective bundle associated to $E$. We
denote the tautological line bundle on $P(E)$ by
$\mathcal{O}_{P(E)}(1)$ and first Chern class by $c_1$ (see
\cite{FG}).
\begin{enumerate}
\item (Friedlander-Gabber) If $p\geq r$, then the map
$$\Psi\equiv \sum^r_{j=0}c_1(O_{P(E)}(1))^{r-j}\circ \pi^*:
\bigoplus^r_{j=0} Z_{p-j}(Y) \rightarrow Z_p(P(E))$$ is a homotopy
equivalent.

\item If $0\leq p<r$, then the map
$$\Psi:
\bigoplus^r_{j=0} Z_{p-j}(Y) \longrightarrow Z_p(P(E))$$ is a
homotopy equivalent where $\Psi\equiv
\sum^p_{j=0}(\pi_1^*)^{-1}\circ c_1(O_{P(E)\times
\C^{r-p}}(1))^{r-j}\circ \pi^* \circ
\pi_2^*+\sum^r_{j=p+1}(\pi_1^*)^{-1}\circ c_1(O_{P(E)\times
\C^{r-p}}(1))^{r-j}\circ \pi^* \circ \pi_3^*$, and $\pi_1:
P(E)\times \C^{r-p} \rightarrow P(E), \pi_2: Y\times \C^{r-p}
\rightarrow Y, \pi_3:Y\times \C^{r-p} \rightarrow Y\times \C^{j-p}$
are the projections.

\item For any $p$, we have a commutative diagram:
$$\xymatrix{\Psi:\bigoplus^r_{j=0}L_{p-j}H_{k-2j}(Y) \ar[r] \ar[d]_{\Phi} &
L_pH_k(P(E))\ar[d]^{\Phi}\\
\Psi':\bigoplus^r_{j=0}H^{BM}_{k-2j}(Y) \ar[r] & H^{BM}_k(P(E))}$$
where $\Psi'$ is the counterpart of $\Psi$ in singular homology.
\end{enumerate}

\label{projective bundle}
\end{proposition}

\begin{proof*}
We prove (ii). We have two homotopy equivalences:
$$\bigoplus^p_{j=0}Z_{p-j}(Y) \overset{\pi_2^*}{\longrightarrow}
\bigoplus^p_{j=0}Z_{r-j}(Y\times \C^{r-p}),$$ and
$$\bigoplus^r_{j=p+1}Z_0(Y\times\C^{j-p})\overset{\pi_3^*}{\longrightarrow}
\bigoplus^r_{j=p+1}Z_{r-j}(Y\times \C^{r-p})$$ Combining them
together, we get
$$\bigoplus^r_{j=0}Z_{p-j}(Y)= \bigoplus^p_{j=0}Z_{p-j}(Y)\oplus\bigoplus^r_{j=p+1}Z_0(Y\times
\C^{j-p}) \cong \bigoplus^r_{j=0}Z_{r-j}(Y\times \C^{r-p})$$

Consider $\C^{r-p}$ as a zero rank vector bundle over itself, then
$E\times \C^{r-p}$ is an algebraic vector bundle over $Y\times
\C^{r-p}$ of rank $r+1$. We have $P(E\times \C^{r-p})=P(E)\times
\C^{r-p}$. By (i), we have a homotopy equivalence
$$\sum^r_{j=0}c_1(O_{P(E)\times \C^{r-p}}(1))^{r-j}\circ \pi^*:\bigoplus^r_{j=0} Z_{r-j}(Y\times \C^{r-p})
\longrightarrow Z_r(P(E)\times \C^{r-p})$$

Combining with the homotopy equivalence
$(\pi_1^*)^{-1}:Z_r(P(E)\times \C^{r-p}) \rightarrow Z_p(P(E))$, we
are done. (iii) follows from the fact that all the maps in (i) and
(ii) are induced from algebraic maps and $\Phi$ is a natural
transformation.
\end{proof*}

The following is the blow-up formula in Lawson homology from
\cite{Hu}. We state the result for integral coefficients but we use
the formula only in rational coefficients.

\begin{proposition}
Let $X$ be a smooth projective variety and $i':Y\hookrightarrow X$ a
smooth subvariety of codimension $r+1$. Let $\sigma: Bl_YX
\rightarrow X$ be the blow-up of $X$ along $Y$,
$\pi:E=\sigma^{-1}(Y)\rightarrow Y$ the projection, and $i:E
\rightarrow Bl_YX$ the inclusion map. For $p\geq 0$ and $k\geq 2p$,
\begin{enumerate}
\item
the map
$$I_{p, k}:\bigoplus^r_{j=1}L_{p-j}H_{k-2j}(Y)\oplus L_pH_k(X)
\rightarrow L_pH_k(Bl_YX)$$ defined by
$$I_{p, k}(u_1,..., u_r,
u)=\sum^r_{j=1}i_*h^{r-j}\pi^*u_j+\sigma^*u$$ is an isomorphism
where $h\in L_{m-1}H_{2(m-1)}(E)$ is the class defined by a
hyperplane section of $E$.

\item
There is a split short exact sequence
$$0\rightarrow L_pH_k(Y)\overset{\psi_1}{\rightarrow} L_pH_k(E)\oplus L_pH_k(X) \overset{\psi_2}{\rightarrow}L_pH_k(Bl_YX)
\rightarrow 0$$ where $\psi_1(x)=(h^r\pi^*x, -i'_*x)$, and
$\psi_2(x)=(\widetilde{x}, y)=i_*\widetilde{x}+\sigma^*y$.
\end{enumerate}
\label{blowup formula}
\end{proposition}

\begin{definition}
Suppose that $X$ is a smooth projective variety. Let
$T_pH_k(X):=\Phi(L_pH_k(X;\Q))\subset H_k(X; \Q)$ be the image of
the natural transformation $\Phi$ from Lawson homology to singular
homology and $T_{p, k}=dim\Phi(L_pH_k(X; \Q))$ be its dimension for
$k\geq 2p\geq 0$.
\end{definition}

Combine with the blow-up formula in singular homology (see
\cite{GH}, Chapter 4.6), and as an immediate consequence of the
Proposition above, we get the following crucial equality.

\begin{proposition}
Suppose that $X$ is a smooth projective variety. Then we have
$$T_{p, k}(Bl_YX)-T_{p, k}(E(Y))=T_{p, k}(X)-T_{p, k}(Y)$$
where $Bl_YX$ is the blow-up of $X$ along a smooth subvariety $Y$
and $E(Y)$ is the exceptional divisor.
\end{proposition}

Since $K_0(SPV)=\Z(SPV)/\sim_{bl}$ where $\sim_{bl}$ is the blow-up
relation, we see that $T_{p, k}$ induces a group homomorphism from
$K_0(SPV)$ to $\Z$.

\begin{proposition}
For a smooth projective variety $X$, define
$$T_{p, k}(X\cdot
\L^i):=T_{p, k}(X\times \P^i)-T_{p, k}(X\times \P^{i-1})$$ for
$k\geq 2p\geq 0$.

\begin{enumerate}
\item
The map $T_{p, k}$ induces a group homomorphism from $\mathcal{M}$
to $\Z$.

\item The map $T_{p, k}$ induces a group homomorphism from $\overline{\mathcal{M}}[(\frac{1}{\L^i-1})_{i\geq
1}]$ to $\Z$.
\end{enumerate}

\label{homomorphism}
\end{proposition}

\begin{proof*}
\begin{enumerate}
\item
From the projective bundle theorem for trivial bundles and applying
the natural transformations from Lawson homology to singular
homology, we have the following commutative diagram:
$$\xymatrix{  L_pH_k(X\times \P^i; \Q)\ar[d]_{\Phi} &\cong &  L_pH_k(X\times \P^{i-1};
\Q) \oplus L_{p-i}H_{k-2i}(X; \Q) \ar[d]^{\Phi} \\
H_k(X\times \P^i; \Q)& \cong & H_k(X\times \P^{i-1}; \Q)  \oplus
H_{k-2i}(X; \Q)}$$

Then $T_{p, k}(X\cdot \L^i)=T_{p, k}(X\times \P^i)-T_{p, k}(X\times
\P^{i-1})= T_{p-i, k-2i}(X)$.

For $\frac{A}{\L^i}\in \mathcal{M}$ where $A\in K_0(SPV)$, we define
$$T_{p, k}(\frac{A}{\L^i}):=T_{p+i, k+2i}(A)$$

Since $T_{p, k}(\frac{A\cdot \L^j}{\L^{i+j}})=T_{p+i+j,
k+2(i+j)}(A\cdot \L^j)=T_{p+i, k+2i}(A)=T_{p, k}(\frac{A}{\L^i})$,
$T_{p, k}$ is well defined over $\mathcal{M}$. Extending $T_{p, k}$
by linearity, we get a group homomorphism from $\mathcal{M}$ to
$\Z$.

\item
The kernel of the canonical map $\phi:\mathcal{M} \rightarrow
\hat{\mathcal{M}}$ is $\cap_n\mathcal{F}^n\mathcal{M}$. If $A\in
\cap_n\mathcal{F}^n\mathcal{M}$, then
$A=\sum_ib_i\frac{B_i}{\L^{n_i}}$ for some $b_i\in \Z, B_i\in
K_0(SPV)$ such that $n_i-dim B_i\geq 1$. Therefore $T_{p,
k}(A)=\sum_ib_iT_{p, k}(\frac{B_i}{\L^{n_i}})= \sum_ib_iT_{p+n_i,
k+2n_i}(B_i)=0$ since $n_i>\mbox{ dim }B_i$. Hence $T_{p, k}$
induces a group homomorphism from
$\overline{\mathcal{M}}=\phi(\mathcal{M})\subset \hat{\mathcal{M}}$
to $\Z$.

Write $\frac{1}{\L^i-1}=\sum^{\infty}_{j=0}a_j\L^j$. For $[X]\in
\overline{\mathcal{M}}$, define
$$T_{p, k}(X(\L^i-1)^{-1}):=\sum^{\infty}_{j=0}a_jT_{p,k}(X\cdot
\L^j)$$ which is a finite sum. Hence $T_{p, k}$ extends to a group
homomorphism from $\overline{\mathcal{M}}[(\frac{1}{\L^i-1})_{i\geq
1}]$ to $\Z$.
\end{enumerate}
\end{proof*}

\begin{proposition}
For a smooth projective variety $X$ of dimension $m$, let $h^{p,
q}(X)$ be the $(p, q)$-Hodge number of $X$ and $h_{m-p, m-q}(X)$ be
the dimension of the Poincar\'e dual of $H^{p, q}(X)$ in the
homology group $H_{2m-(p+q)}(X; \C)$. Then $h_{p, q}$ and $h^{p, q}$
induce group homomorphisms from
$\overline{\mathcal{M}}[(\frac{1}{\L^i-1})_{i\geq 1}]$ to $\Z$.
\label{hodge numbers}
\end{proposition}

\begin{proof*}
We show this for $h_{p, q}$. The maps in the exact sequence
$$0\rightarrow H_n(X; \C) \rightarrow H_n(Bl_YX; \C)\oplus H_n(Y; \C)
\rightarrow H_n(E(Y); \C) \rightarrow 0$$ are easy to see to be
morphisms of Hodge structures, hence we have
$$h_{p, q}(Bl_YX)-h_{p, q}(E(Y))=h_{p, q}(X)-h_{p, q}(Y)$$
which implies that $h_{p, q}$ induces a group homomorphism from
$K_0(SPV)$ to $\Z$.

Define
$$h_{p, q}(X\cdot \L^i):=h_{p, q}(X\times \P^i)-h_{p,
q}(X\times \P^{i-1}).$$

From the isomorphism $H_{p, q}(X\times \P^i)\cong H_{p, q}(X\times
\P^{i-1})\oplus H_{p-i, q-i}(X)$, we get
$$h_{p, q}(X\cdot \L^i)=h_{p-i, q-i}(X)$$
and then we extend $h_{p, q}$ as in the Proposition above.
\end{proof*}

We recall that the niveau filtration $\{N_pH_*(X; \Q)\}_{p\geq 0}$
of $H_*(X; \Q)$ is defined by
$$N_pH_k(X; \Q)=\mbox{ span }\{ \mbox{ images } i_*: H_k(Y; \Q) \rightarrow H_k(X;
\Q)|i:Y \hookrightarrow X, dimY \leq p\}$$ The geometric filtration
$\{G_pH_*(X; \Q)\}_{p\geq 0}$ of $H_*(X; \Q)$ is defined by
$$G_pH_k(X; \Q)=N_{k-p}H_k(X; \Q)$$
We define the homological Hodge filtration to be
$$F_pH_k(X; \C):=\bigoplus_{t\leq p}H_{t, k-t}(X)$$
and define the homological rational Hodge filtration to be
$$F^h_pH_k(X; \Q)=\mbox{ largest sub-Hodge structure of }
F_pH_k(X; \C)\cap H_k(X; \Q)$$ The homological generalized Hodge
conjecture says that for a smooth projective variety $X$,
$$F^h_pH_k(X; \Q)=N_pH_k(X; \Q)$$
The Friedlander-Mazur conjecture (see \cite{FM}) says that
$$T_pH_k(X; \Q)=G_pH_k(X; \Q)$$ and the Friedlander-Lawson
conjecture says that
$$L_pH_k(X; \Q)\rightarrow H_k(X; \Q)$$ is surjective if $k\geq m+p$
where $m$ is the dimension of $X$. This conjecture was proved by the
author in \cite{T} by assuming the Grothendieck standard conjecture.

\begin{definition}
We say that a statement is a $K$-statement if it is true for a
smooth projective variety $X$, then it is true for all varieties
which are $K$-equivalent to $X$.
\end{definition}

We will show that all these conjectures are $K$-statements.

\begin{proposition}
Let $G_{j, n}(X):=dim G_jH_n(X; \Q), F_{j, n}(X):= dim F^h_jH_n(X;
\Q)$. Then $G_{j, n}$ and $F_{j, n}$ extend to
$\overline{\mathcal{M}}[(\frac{1}{\L^i-1})_{i\geq 1}]$.
\end{proposition}

\begin{proof*}
By Lemma 2.3 of \cite{AK}, we have a short exact sequence of pure
Hodge structures:
$$0\rightarrow H_n(X; \C) \rightarrow H_n(Bl_YX; \C)\oplus H_n(Y; \C)
\rightarrow H_n(E(Y); \C) \rightarrow 0$$ which give us the
following formula
$$F_{j, n}(Bl_YX)+F_{j, n}(Y)=F_{j, n}(X)+F_{j, n}(E(Y)).$$
We define
$$F_{j, n}(X\cdot \L^k):=F_{j, n}(X\times \P^k)-F_{j,
n}(X\times \P^{k-1})$$ Since $X\times \C^k=X\times \P^k- X\times
\P^{k-1}$, from the mixed Hodge theory, there is a long exact
sequence of mixed Hodge structures:
$$ \cdots \rightarrow H_n(X\times \P^{k-1}; \C) \rightarrow H_n(X\times
\P^k; \C) \rightarrow H^{BM}_n(X\times \C^k; \C) \rightarrow
H_{n-1}(X\times \P^{k-1}; \C) \rightarrow \cdots$$ but the map
induced by inclusion $H_i(X\times \P^{k-1}; \C)\rightarrow
H_i(X\times \P^k; \C)$ is always an injection, and therefore we get
an exact sequence:
$$0\rightarrow H_n(X\times \P^{k-1}; \C) \rightarrow H_n(X\times \P^k;
\C) \rightarrow H^{BM}_n(X\times \C^k; \C)\rightarrow 0$$ The
isomorphism $H_{n-2k}(X; \C) \cong H^{BM}_n(X\times \C^k; \C)$ is an
isomorphism of Hodge structures of type $(k, k)$, hence
$$F_{j, n}(X\cdot \L^k)=\mbox{ dim }F^h_jH^{BM}_n(X\times \C^k; \Q)=
\mbox{ dim }F^h_{j-k}H_{n-2k}(X; \C) =F_{j-k, n-2k}(X)$$ Then extend
as in the Proposition \ref{homomorphism}, we get a group
homomorphism from $\overline{\mathcal{M}}[(\frac{1}{\L^i-1})_{i\geq
1}]$ to $\Z$.

By a homological version of Lemma 2.4 of \cite{AK}, we get
$$G_{j, n}(Bl_YX)+G_{j, n}(Y)=G_{j, n}(X)+G_{j, n}(E(Y)).$$
Define $$G_{j, n}(X\cdot \L^k):=G_{j, n}(X\times \P^k)-G_{j,
n}(X\times \P^{k-1}).$$ Since a morphism of Hodge structures
preserve niveau filtration, we have $G_{j, n}(X\cdot \L^k)=G_{j-k,
n-2k}(X)$. Then similar to the construction above, $G_{j, n}$
extends to a group homomorphism from
$\overline{\mathcal{M}}[(\frac{1}{\L^i-1})_{i\geq 1}]$ to $\Z$.
\end{proof*}

By composing with the isomorphism in Proposition 1\ref{isomorphism},
we have the following crucial result.
\begin{theorem}
The group homomorphisms $T_{p, k}, h^{p, k}, h_{p, k}, G_{p, k},
F_{p, k}$ induce group homomorphisms from
$\overline{\mathcal{N}}[(\frac{1}{\L^i-1})_{i\geq 0}]$ to $\Z$. We
will abusively use the same notations for the induced homomorphisms.
\end{theorem}

\begin{remark}
For a quasi-projective variety $X$, we may define $T'_{p, k}(X)$ to
be the dimension of $\Phi(L_pH_k(X; \Q))\subset H_k^{BM}(X; \Q)$.
Even though this definition makes sense, in general it does not
equal to $T_{p, k}(X)$. For example $T_{0, 0}(\C^*)=T_{0,
0}(\P^1)-T_{0, 0}(\{0, \infty\})=1-2=-1$ but $T'_{0, 0}(\C^*)$ is a
nonnegative integer.
\end{remark}

We recall that the Lawson homology group $L_pH_{2p}(X)=$ algebraic
$p$-cycles quotient by algebraic equivalence and the natural
transformation $\Phi:L_pH_{2p}(X)\rightarrow H_{2p}(X)$ is the cycle
map. Hence $\Phi(L_pH_{2p}(X))$ is the subspace of $H_{2p}(X)$
generated by algebraic cycles. We recall that the Grothendieck
standard conjecture A (GSCA) predicts that for a smooth projective
variety $X$ of dimension $m$, $T_{p, 2p}(X)=T_{m-p, 2(m-p)}(X)$
where $p\leq [\frac{m}{2}]$.

\begin{proposition}
If $Bl_YX$ is the blow-up of a smooth projective variety $X$ at a
smooth center $Y$ of codimension $r+1$ and if the GSCA is true for
$Y$, then the GSCA on $Bl_YX$ is equivalent to the GSCA on $X$.
\end{proposition}

\begin{proof*}
Let the dimension of $X$ be $m$ and $p\leq [\frac{m}{2}]$. Then the
dimension of $Y$ is $m-r-1$. Let $A_p(Y)=\sum^r_{j=1}T_{p-j,
2(p-j)}(Y)$, $B_p(X)=T_{p, 2p}(X)$, $B_p(Bl_YX)=T_{p, 2p}(Bl_YX)$.
We have $A_p(Y)+B_p(X)=B_p(Bl_YX)$.

From the calculation $A_{m-p}(Y)= \sum^r_{j=1}T_{m-p-j,
2(m-p-j)}(Y)= \sum^r_{j=1}T_{p+j-r-1, 2(p+j-r-1)}(Y)=$ \\
$\sum^r_{j=1}T_{p-j, 2(p-j)}(Y)=A_p(Y)$, we get
$B_p(X)-B_{m-p}(X)=B_p(Bl_YX)-B_{m-p}(Bl_YX)$ which means that the
GSCA holds on $X$ if and only if it holds on $Bl_YX$.
\end{proof*}

\begin{corollary}
If the GSCA holds for smooth projective varieties of dimension less
than $m-1$, then the GSCA is a birational statement for smooth
projective varieties of dimension $m$.
\end{corollary}

\begin{proof*}
By the Weak Factorization Theorem of birational maps (see
\cite{AKMW}), we are able to decompose a proper birational map as a
sequence of blowing-ups and blowing-downs, then we apply the result
above.
\end{proof*}

Since we know that the GSCA is true for smooth varieties of
dimension less than or equal to 4, we have the following result.
\begin{corollary}
The GSCA is invariant under birational equivalence of smooth
varieties of dimension less than 7.
\end{corollary}

For a projective manifold $X$, let $N_{j, n}(X):=\mbox{ dim
}N_jH_n(X; \Q)$.

\begin{proposition}
If $Bl_YX$ is the blow-up of a smooth projective variety $X$ at a
smooth center $Y$ of codimension $r+1$ and the generalized Hodge
conjecture is true for $Y$, then the generalized Hodge conjecture on
$Bl_YX$ is equivalent to the generalized Hodge conjecture on $X$.
\end{proposition}

\begin{proof*}
We have $N_{j, n}(Bl_YX)=N_{j, n}(X)+\sum^r_{i=1}N_{j-i, n-2i}(Y)$
and $F_{j, n}(Bl_YX)=F_{j, n}(X)+\sum^r_{i=1}F_{j-i, n-2i}(Y)$. By
the assumption that the generalized Hodge conjecture is true for
$Y$, we have $N_{j, n}(Bl_YX)-F_{j, n}(Bl_YX)=N_{j, n}(X)-F_{j,
n}(X)$. This completes the proof.
\end{proof*}

Again by using the Weak Factorization Theorem of birational maps, we
get the following result. We do not know who is the first to have
this result, but a proof without using the Weak Factorization
Theorem can be found in \cite{Ara1}.
\begin{corollary}
If the Hodge conjecture is true for dimension less than $m-1$, then
the Hodge conjecture is a birational statement for smooth varieties
of dimension $m$. In particular it is a birational statement for
dimension less than 6.
\end{corollary}

Since two $K$-equivalent varieties have the same image in
$\overline{\mathcal{N}}$, any group homomorphism defined previously
gives the same value at them. Then the following result is an
immediate consequence.

\begin{theorem}
The Friedlander-Mazur conjecture, the Friedlander-Lawson conjecture,
the Grothendieck standard conjecture and the generalized Hodge
conjecture are $K$-statements. \label{conjectures}
\end{theorem}

For the case of generalized Hodge conjecture, this result was proved
by Arapura and Kang in \cite{AK}. By a result of Wang (see
\cite{Wang}, Corollary 1.10), two birational smooth minimal models
are $K$-equivalent, hence in particular we have the following
result.

\begin{corollary}
If any conjecture in Theorem \ref{conjectures} is true for a smooth
minimal model, then it is true for any smooth minimal model which is
birational to it.
\end{corollary}

\section{Stringy functions}
\begin{definition}
A motivic invariant is a group homomorphism from
$\overline{\mathcal{N}}[(\frac{1}{\L^i-1})_{i\geq 1}]$ to $\Z$.
\end{definition}

We have seen several motivic invariants: $T_{j, n}$, $G_{j, n}$,
$F_{j, n}$, $h^{p, q}$ and $h_{p, q}$. One of the most important
properties of these invariants is that they satisfy $\phi_{j,
n}(X\times \L^k)=\phi_{j-ak, n-bk}(X)$ for some numbers $a, b$. This
enables us to associate a stringy $\phi$-function to $\phi$. Before
we consider the general case, let us exemplify this by Batyrev's
stringy $E$-function.

\begin{example}
Let us recall some definitions from \cite{Baty}. For a variety $X$
of dimension $m$, let
$$e^{p, q}(X):=\sum_{0\leq k\leq 2m}(-1)^kh^{p, q}(H^k_c(X;
\C))$$ where $h^{p, q}(H^k_c(X; \C))$ is the $(p, q)$-Deligne-Hodge
number of the cohomology groups with compact support of $X$. For a
projective manifold, the number $e^{p, q}$ is same as the Hodge
number $h^{p, q}$. The $E$-polynomial $E(X; u, v)\in \Z[u, v]$ is
defined to be
$$E(X; u, v):=\sum_{p, q}e^{p, q}(X)u^pv^q.$$
This is a finite sum and $E(X\times \L^k; u, v)=(uv)^kE(X; u, v)$.
Therefore by defining $E(\L^{-1}; u, v)=(uv)^{-1}$, we are able to
extend $E$ to a group homomorphism
$E:\overline{\mathcal{N}}[(\frac{1}{\L^i-1})_{i\geq 1}] \rightarrow
\Z[u, v, (uv)^{-1}]$.

If $X$ is a normal irreducible algebraic variety with at worst
log-terminal singularities, and $\rho: Y \rightarrow X$ is a
resolution of singularities such that the relative canonical divisor
$D=\sum^r_{i=1}a_iD_i$ has simple normal crossings. Then the stringy
$E$-function of $X$ is defined to be
$$E_{st}(X; u, v):=\sum_{J\subset I}E(D^0_J; u, v)\prod_{j\in
J}\frac{uv-1}{(uv)^{a_j+1}-1}$$ where $I=\{1,..., r\}$. If
$E_{st}(X; u, v)=\sum_{p, q}a_{p, q}u^pv^q$ is a polynomial, we
define the stringy Hodge numbers of $X$ to be
$$h^{p, q}_{st}(X):=(-1)^{p+q}a_{p, q}.$$
\end{example}

Now we come to the general case.

\begin{definition}
We say that a family of motivic invariants $\phi=\{\phi_{j, n}|j,
n\in \Z\}$ is of type $(a, b)\in \Z\times \Z$ if $\phi_{j,
n}(X\times \L^k)=\phi_{j-ak, n-bk}(X)$ for any $j, n$ and any
varieties $X$. And we say that $\phi$ is bounded if $\phi_{j, n}(X)$
vanishes for $j, n$ large enough, depending on $X$.
\end{definition}

For example $T=\{T_{j, n}|j, n\in \Z\}$ is of type $(1, 2)$ and
$h=\{h^{p, q}|p, q\in \Z\}$ is of type $(1, 1)$ where $T_{j, n},
h^{p, q}$ are defined to be zero if any $j, n, p, q$ is negative.

\begin{definition}
Suppose that $\phi=\{\phi_{j, n}|j, n\in \Z\}$ is a family of
bounded motivic invariants of type $(a, b)$, then define
$$\phi(X; u, v):=\sum_{j, n}\phi_{j, n}(X)u^jv^n$$ and $$\phi(\L^{-1};
u, v):=(u^av^b)^{-1},$$ we get a group homomorphism
$$\phi:\overline{\mathcal{N}}[(\frac{1}{\L^i-1})_{i\geq 1}] \rightarrow \Z[u, v,
(u^av^b)^{-1}].$$

If $X$ is a normal irreducible algebraic variety with at worst
log-terminal singularities, and $\rho: Y \rightarrow X$ is a
resolution of singularities such that the relative canonical divisor
$D=\sum^r_{i=1}a_iD_i$ has simple normal crossings. Then the stringy
$\phi$-function associated to $\phi$ is defined to be
$$\phi^{st}(X; u, v):=\sum_{J\subset I}\phi(D^0_J; u, v)\prod_{j\in
J}\frac{u^av^b-1}{(u^av^b)^{a_j+1}-1}$$ where $I=\{1,..., r\}$. If
$X$ is projective and $\phi^{st}(X; u, v)=\sum_{p, q}a_{p, q}u^pv^q$
is a polynomial, we define the $(p, q)$-stringy $\phi$-numbers of
$X$ to be
$$\phi^{st}_{p, q}(X):=(-1)^{p+q}a_{p, q}$$
\end{definition}

\begin{proposition}
The stringy $\phi$-numbers of $X$ defined above are independent of
resolution of singularities.
\end{proposition}

\begin{proof*}
Let $\rho_1:X_1 \rightarrow X, \rho_2:X_2 \rightarrow X$ be two
resolution of singularities. Take another resolution of
singularities $\alpha:Y \rightarrow X$ which dominates $\rho_1,
\rho_2$, i.e., we have the following commutative diagram:
$$\xymatrix{& \ar[ld]_{\alpha_1} Y \ar[dd]^{\alpha} \ar[rd]^{\alpha_2} &\\
X_1 \ar[rd]_{\rho_1} & & X_2 \ar[ld]^{\rho_2}\\
& X &}$$ Let $K_{X_1}=\rho^*_1K_X+D_1, K_{X_2}=\rho^*_2K_X+D_2$ and
$D=K_Y-\alpha^*K_X$. Then
$\alpha^*_1D_1+K_{Y|X_1}=\alpha^*_1D_1+K_Y-\alpha^*_1K_{X_1}=K_Y-\alpha^*K_X=\alpha^*_2D_2+K_{Y|X_2}$.
Therefore by the change of variables formula,
$$\int_{\mathcal{L}(X_1)}\L^{-ord_tD_1}=\int_{\mathcal{L}(Y)}\L^{-ord_t(\alpha^*_1D_1+K_{Y|X_1})}=
\int_{\mathcal{L}(Y)}\L^{-ord_t(\alpha^*_2D_2+K_{Y|X_2})}=\int_{\mathcal{L}(X_2)}\L^{-ord_tD_2}$$
Taking $\phi^{st}$ on both sides, this shows that $\phi^{st}$ is
independent of resolution of singularities.
\end{proof*}

With all these definitions, we are able to ask the stringy version
of some conjectures.

\begin{conjecture}
Suppose that $X$ is a $m$-dimensional normal irreducible projective
variety with at worst log-terminal singularities. Let $T^{st}_p,
T^{st}_{p, q}, G^{st}_{p, q}, F^{st}_{p, q}$ be the $(p, q)$-stringy
numbers of the families $T=\{T_{p, 2p}|p\in \Z\}, T'=\{T_{p, q}|p,
q\in \Z\}, G=\{G_{p, q}|p, q\in \Z\}, F=\{F_{p, q}|p, q\in \Z\}$
respectively.
\begin{enumerate}
\item (Stringy GSCA) Is $T^{st}_{p, 2p}(X)=T^{st}_{m-p, 2(m-p)}(X)$?
\item (Stringy morphic conjecture) Is $T^{st}_{p, q}(X)=T^{st}_{m-p,
2m-q}(X)$?
\item (Stringy generalized Hodge conjecture) Is $G^{st}_{p,
q}(X)=F^{st}_{p, q}(X)$?
\item (Stringy Hodge conjecture) Is $T^{st}_{p, 2p}(X)=F^{st}_{p, 2p}(X)$?
\end{enumerate}
\end{conjecture}

By \cite[Theorem 3.7]{Baty}, for a projective $\Q$-Gorenstein
variety $X$ of dimension $d$ with at worst log-terminal
singularities, Batyrev's stringy $E$-function satisfies the
equality:
$$E_{st}(X; u, v)=(uv)^dE_{st}(X; u^{-1}, v^{-1})$$
this follows basically from the strong Lefschetz theorem. Similar
calculation shows that $G^{st}(X; u, v)=(uv^2)^dG^{st}(X; u^{-1},
v^{-1})$ and
$$F^{st}(X; u, v)=(uv^2)^dF^{st}(X; u^{-1}, v^{-1})$$ which follows
from the fact that the Lefschetz isomorphism is an isomorphism of
Hodge structures. These facts are some special cases of the
following conjecture.

\begin{conjecture}(Generalized stringy GSCA)
Let $X$ be as above. If $\phi$ is a family of bounded motivic
invariants of type $(a, b)$, then $\phi^{st}(X; u,
v)=(u^av^b)^d\phi^{st}(X; u^{-1}, v^{-1})$.
\end{conjecture}

If $X$ is smooth projective, this is just the GSCA.

We verify this conjecture for normal projective $\Q$-Gorenstein
toric varieties.

\begin{theorem}
Suppose that $X$ is a normal projective $\Q$-Gorenstein toric
varieties, then the generalized stringy GSCA holds.
\end{theorem}

\begin{proof*}
Let $d$ be the dimension of $X$ and $\phi=\{\phi_{i, j}\}$ be a
family of bounded motivic invariants of type $(a, b)$. By
\cite[Theorem 3.7]{Baty}, the stringy $E$-function of $X$ satisfies
the following relation: $E_{st}(X; u, v)=(uv)^dE_{st}(X; u^{-1},
v^{-1})$. And by \cite[Theorem 4.3]{Baty}, $E_{st}(X; u,
v)=\sum_{\sigma\in \sum}\sum_{n\in \sigma^0\cap
\N}(uv)^{-\varphi(n)}$ where $X$ is defined by the fan $\sum$ on the
lattice $\N$, and $\varphi$ is a supporting function of $X$.
Comparing the equality of the $E$-function, we get
$(-1)^d\sum_{\sigma\in \sum}\sum_{n\in \sigma^0\cap
N}(uv)^{\varphi(n)}=\sum_{\sigma\in \sum}\sum_{n\in \sigma^0\cap
\N}(uv)^{-\varphi(n)}$. Now follow exactly the same calculation as
in \cite[Theorem 4.3]{Baty}, the stringy function satisfies the
equality: $\phi^{st}(X; u, v)=(u^av^b-1)^d\sum_{\sigma\in
\sum}\sum_{n\in \sigma^0\cap
\N}(uv)^{-\varphi(n)}=(-1)^d(u^av^b-1)^d\sum_{\sigma\in
\sum}\sum_{n\in \sigma^0\cap
\N}(uv)^{\varphi(n)}=(u^av^b)^d\phi^{st}(X; u^{-1}, v^{-1})$.
\end{proof*}

\begin{proposition}
For a toric variety $X_{N, \Sigma}$ of dimension $m$,
\begin{enumerate}
\item $$h_{p, q}(X_{N, \Sigma})=
\left\{
       \begin{array}{ll}
       \sum^m_{k=0}d_{m-k}(-1)^{k-p}\binom{k}{k-p}, & \hbox{ if } p=q \hbox{ and } 0\leq p \leq m, \\
       0, & \hbox{ otherwise.}
       \end{array}
\right.
$$
where $d_k$ is the number of cones of dimension $k$ in $\Sigma$.
\item the number
$e^{p, q}(X_{N, \Sigma})$ is equal to $h_{p, q}(X_{N, \Sigma})$.

\item $T_{p, 2p}(X_{N, \Sigma})=T_{p, k}(X_{N, \sigma})=N_{p, k}(X_{N, \Sigma})=F_{p,
k}(X_{N, \Sigma})=h_{p, p}(X_{N, \Sigma})$ for $k\geq 2p$. In
particular, the Friedlander-Lawson conjecture, the Friedlander-Mazur
conjecture, the generalized Hodge conjecture, the Grothendieck
standard conjecture are true for smooth toric varieties.

\end{enumerate}
\end{proposition}

\begin{proof*}
\begin{enumerate}
\item
The action of the torus $\mathbb{T}\sim (\C^*)^m$ on $X_{N, \Sigma}$
induces a stratification of $X_{N, \Sigma}$ into orbits of the torus
action $O_{\tau}\cong (\C^*)^{m-dim\tau}$, one for each cone
$\tau\in \Sigma$. Then we have
$$[X_{N, \Sigma}]=\sum^m_{k=0}d_{m-k}[\L-1]^k$$ in $K_0(SPV)$ . Since
$[\L-1]^k=\sum^k_{i=0}\binom{k}{i}(-1)^k\L^{k-i}$, we have
$$h_{p, q}([\L-1]^k)=\left\{
             \begin{array}{ll}
               \binom{k}{k-p}(-1)^{k-p}, & \hbox{ if } p=q \hbox{ and } 0\leq p \leq k \\
               0, & \hbox{ otherwise}
             \end{array}
           \right.
$$
and substitute into the formula $$h_{p, q}([X_{N,
\Sigma}])=\sum^m_{k=0}d_{m-k}h_{p, q}([\L-1]^k),$$ then we get the
result.

\item It was calculated by Batyrev (see \cite{Baty}) that $E(X_{N,
\Sigma}; u, v)=\sum^m_{k=0}d_{m-k}(uv-1)^k$. Then we make a simple
comparison to the coefficients of $E(X_{N, \sigma})$ with the
corresponding $h_{p, q}(X_{N, \Sigma})$.

\item We note that $T_{p, 2p}(\L^k)=h_{p, p}(\L^k)$ for any $p, k$,
this implies the equality of $T_{p, 2p}(X_{N, \Sigma})=h_{p,
p}(X_{N, \Sigma})$. The number $h_{p, q}(X_{N, \Sigma})=0$ if $p\neq
q$, this implies that $T_{p, 2p}(X_{N, \Sigma})=T_{p, k}(X_{N,
\sigma})=N_{p, k}(X_{N, \Sigma})=F_{p, k}(X_{N, \Sigma})=h_{p,
p}(X_{N, \Sigma})$ for $k\geq 2p$. For smooth toric varieties,
$T_{p, 2p}(X_{N, \Sigma})=h_{p, p}(X_{N, \Sigma})$ means that the
homology group $H_{2p}(X_{N, \Sigma}; \Q)$ is generated by algebraic
cycles hence all the conjectures are trivially true.
\end{enumerate}
\end{proof*}

In their paper \cite{BB}, Batyrev and Borisov proved the mirror
duality conjecture for stringy Hodge numbers of Calabi-Yau complete
intersections in Gorenstein Fano toric varieties, i.e., for a mirror
pair $(V, W)$ of such varieties of dimension $n$, their stringy
$E$-functions satisfies the relation
$$E_{st}(V; u, v)=(-u)^nE_{st}(W; u^{-1}, v)$$ which in particular
gives the rotation of the Hodge diamond: $h^{p, q}(V)=h^{n-p,
q}(W)$. We wonder if similar relation is true for a family of
bounded motivic invariants. We form our conjecture below.

\begin{conjecture}
Given a family of bounded motivic invariants $\phi=\{\phi_{j, n}|j,
n\in \Z\}$ of type $(a, b)$. Then if $(V, W)$ is a mirror pair from
Batyrev-Borisov's construction, then
$$\phi^{st}(V; u, v)=(-u^a)^n\phi^{st}(W; u^{-1}, v).$$
\end{conjecture}

\section{Lawson-Deligne-Hodge polynomials}

\subsection{Varieties with finitely generated Lawson homology
groups} Let $VarFL$ be the collection of all quasi-projective
varieties $X$ such that the dimension of $L_rH_n(X; \Q)$ is finite
for all nonnegative integers $n, r$. Let $X, Y\in VarFL$ and $Y$ be
a locally closed subvariety of $X$. From the localization sequence
of Lawson homology,
$$\cdots \rightarrow L_rH_{n+1}(X; \Q) \rightarrow L_rH_{n+1}(X-Y; \Q) \rightarrow
L_rH_n(Y; \Q) \rightarrow L_rH_n(X; \Q) \rightarrow \cdots$$ we see
that $X-Y$ is also in $VarFL$. Hence we may form the Grothendieck
group $K_0(VarFL)$ of $VarFL$. The ring $\mathscr{L}=\Z\{\L^i|i\in
\Z_{\geq 0}\}$ acts on $K_0(VarFL)$ and we consider $K_0(VarFL)$ as
a $\mathscr{L}$-module under this action. Let $S=\{\L^i\}_{i\geq
0}\subset \mathscr{L}$ and $FL\mathcal{N}=S^{-1}K_0(VarFL)$ be the
localization of $K_0(VarFL)$ with respect to the multiplicative set
$S$. Let $F^kFL\mathcal{N}$ be the subgroup of $FL\mathcal{N}$
generated by elements of the form $\frac{[X]}{\L^i}$ where $i-dim
X\geq k$. Then we have a decreasing filtration
$$\cdots \supset F^kFL\mathcal{N} \supset F^{k+1}FL\mathcal{N}
\supset \cdots$$ of $FL\mathcal{N}$.

\begin{definition}
We definite the $FL$-Kontsevich group to be the completion
$$FL\hat{\mathcal{N}}:=\lim_{\leftarrow
}\frac{\mathcal{N}}{F^kFL\mathcal{N}}$$ with respect to the
filtration defined above.
\end{definition}

Let $X$ be a smooth projective variety of dimension $n$ and $D$ an
effective divisor on $X$ with simple normal crossings. We use the
notation $D^0_J$ as defined in \ref{formula}.

\begin{definition}
We say that a subset $A\subset \mathcal{L}(X)$ is $FL$-cylindrical
if $A=\pi_k^{-1}(C)$ for some $C\subset \mathcal{L}_k(X), C\in
VarFL$. For such set $A$, define
$$\widetilde{\mu}(A):=[C]\L^{-kn}.$$

Let $\mathcal{C}$ be the collection of all countable disjoint unions
of $FL$-cylindrical sets $\coprod_{i\in \N}A_i$ for which
$\widetilde{\mu}(A_i) \rightarrow 0$ in $FL\hat{N}$. An element in
$\mathcal{C}$ is called a $FL$-measurable set. We definite the
$FL$-motivic measure to be $\mu:\mathcal{C}\rightarrow FL\hat{N}$ by
$$\mu(\coprod_{i\in \N}A_i)=\sum_{i\in \N}\widetilde{\mu}(A_i)$$
in $FL\hat{N}$. A function $\alpha:\mathcal{L}(X) \rightarrow \Z\cup
\{\infty\}$ is $FL$-integrable if $\alpha^{-1}(n)$ is
$FL$-measurable for each $n\in \Z\cup \{\infty\}$.
\end{definition}

\begin{proposition}
Let $X$ be a smooth projective variety of dimension $n$ and
$D=\sum^r_{i=1}a_iD_i$ be an effective divisor with simple normal
crossings on $X$. If all $D^0_J$ are in $VarFL$, then $\L^{-ord_t
D}$ is FL-integrable, and $\mu((\L^{-ord_t D})^{-1}(\infty))=0$. We
define the motivic integral of the pair $(X, D)$ to be
$$\int_{\mathcal{L}(X)}\L^{-ord_t D}d\mu:=\sum_{s \in \Z_{\geq
0}}\mu((\mathcal{L}^{-ord_t D})^{-1}(s))$$ in $FL\hat{\mathcal{N}}$.
\label{calculation}
\end{proposition}

This calculation is same as in the proof of \cite[Lemma 2.13]{Craw}
and \cite[Theorem 2.15]{Craw}.

The following result gives a simpler way to see if $D^0_J\in VarFL$.
\begin{proposition}
$D^0_J\in VarFL$ for all $J\subset \{1,..., r\}$ if and only if
$D_J\in VarFL$ for all $J\subset \{1,..., r\}$.
\end{proposition}

\begin{proof*}
Let $R=\{1,..., r\}$. We prove by induction on the length of subsets
of $R$. $D_R=D^0_R\in VarFL$. We assume that for $J\subset R$ with
$|J|>k$, $D_J\in VarFL$. If $|J|=k$, we have
$$D^0_J=D_J-\cup_{i\notin J}D_i=D_J-\cup_{i\notin J}D_i\cap D_J=D_J-\bigcup_{\substack{I\supset J\\|I|=k+1}}D_I$$
Note that for two algebraic varieties $A, B\in VarFL$, if $A\cap
B\in VarFL$, then $A\cup B\in VarFL$. Hence it suffices to prove
that $\bigcup_{\substack{I\supset J\\|I|=k+1}}D_I\in VarFL$. This
follows once we claim that for any subsets $I_1,..., I_t$ of $R$
where each $|I_i|>k$, the union $\cup^t_{i=1}D_{I_i}\in VarFL$. We
use induction again to prove this statement. When $t=1$, $D_{I_1}\in
VarFL$ is by the hypothesis of the first induction. We assume that
this statement is true for $t=n$. Then for $t=n+1$, we have
$$(\bigcup^n_{i=1}D_{I_i})\cap
D_{I_{n+1}}=\bigcup^n_{i=1}(D_{I_i}\cap
D_{I_{n+1}})=\bigcup^n_{i=1}D_{I_i\cup I_{n+1}}\in VarFL$$ By
induction hypothesis, $\cup^n_{i=1}D_{I_i}, D_{I_{n+1}}$ are in
$VarFL$, hence $\cup^{n+1}_{i=1}D_{I_i}\in VarFL$. For another
direction, a similar argument works.
\end{proof*}

We recall a definition from \cite{DL2}.
\begin{definition}
Let $X, Y$ and $F$ be algebraic varieties, and $A\subset X, B\subset
Y$ be constructible subsets of $X$ and $Y$ respectively. We say that
a map $p:A \rightarrow B$ is a piecewise trivial fibration with
fiber $F$, if there exists a finite partition of $B$ in subsets $S$
which are locally closed in $Y$ such that $p^{-1}(S)$ is locally
closed in $X$ and isomorphic, as a variety, to $S\times F$, with $p$
corresponding under the isomorphism to the projection $S\times
F\rightarrow S$.
\end{definition}

By the homotopy property of Lawson homology, we have $L_tH_n(X\times
\C^k)=L_{t-k}H_{n-2k}(X)$ which implies the following Lemma.

\begin{lemma}
For a trivial bundle $X\times \C^k$ over $X$, $X\times \C^k \in
VarFL$ if and only $X\in VarFL$.
\end{lemma}

\begin{proposition}
Let $X, Y$ and $F$ be algebraic varieties, and $A\subset X, B\subset
Y$ be constructible subsets of $X$ and $Y$ respectively. If $p:A
\rightarrow B$ is a piecewise trivial fibration with fibre $\C^k$,
then $A\in VarFL$ if and only if $B\in VarFL$. \label{piecewise
fibration}
\end{proposition}

\begin{proof*}
Consider $B=B_1\coprod B_2$. Then $A=p^{-1}(B_1)\coprod p^{-1}(B_2)$
where $p^{-1}(B_1)\cong B_1\times \C^k, p^{-1}(B_2)\cong B_2\times
\C^k$. From the localization sequences
$$\xymatrix{\cdots \ar[r] & L_pH_n(p^{-1}(B_1)) \ar[r] \ar[d]_{\cong} & L_pH_n(A)
\ar[r] \ar[d] & L_pH_n(p^{-1}(B_2)) \ar[r] \ar[d]_{\cong} & \cdots\\
\cdots \ar[r] & L_{p-k}H_{n-2k}(B_1) \ar[r] & L_{p-k}H_{n-2k}(B)
\ar[r] & L_{p-k}H_{n-2k}(B_2) \ar[r] & \cdots}$$ we see that
$L_pH_n(A)\cong L_{p-k}H_{n-2k}(B)$. Hence $A\in VarFL$ if and only
if $B\in VarFL$. The general case follows by an induction on the
number of components of the partition of $B$.
\end{proof*}

\begin{theorem}(The change of variables formula)
Suppose that $X, Y$ are smooth projective varieties of dimension $d$
and $h:Y \rightarrow X$ is a birational morphism with effective
relative canonical divisor $D=\sum^r_{j=1}a_jD_j$ which has simple
normal crossings. Assume that $X, Y, D^0_J$ are in $VarFL$ for any
$J\subset \{1,..., r\}$, then
$$[X]=\int_{\mathcal{L}(Y)}\L^{-ord_t D}d\mu$$
\label{FL formula}
\end{theorem}

\begin{proof*}
Let $C_k=(\L^{-ord_t D})^{-1}(k)$ for $k\in \Z$.

We claim that $\mu(h_{\infty}(C_k))\in FL\hat{\mathcal{N}}$ and
$$\mu(C_k)=\mu(h_{\infty}(C_k))\cdot \L^k.$$
We have a commutative diagram:
$$\xymatrix{C_k \subset \mathcal{L}(Y) \ar[r]^-{h_{\infty}} \ar[d]_{\pi_t}&
h_{\infty}(C_k)\subset \mathcal{L}(X) \ar[d]^{\pi_t}\\
B'_t\subset \mathcal{L}_t(Y) \ar[r]^{h_t} & B_t\subset
\mathcal{L}_t(X)}$$ where $B'_t=\pi_t(C_k),
B_t=\pi_t(h_{\infty}(C_k))$ are constructible sets. By the
calculation in Proposition \ref{calculation}, we see that
$$B'_t\cong \coprod_{J \subset \{1, ...,
r\}}\coprod_{(m_1, ..., m_r)\in M_{J, k}}D^0_J\cdot\L^{tn-\sum_{j\in
J}m_j}\cdot (\L-1)^{|J|}$$ and by a local calculation in \cite{DL2},
Lemma 3.4(b), the restriction of $h_t$ to $B'_t$ is a piecewise
trivial fibration with fiber $\C^k$ over $B_t$, and from this
calculation we see that
$$B_t\cong \coprod_{J \subset \{1, ...,
r\}}\coprod_{(m_1, ..., m_r)\in M_{J, k}}D^0_J\cdot\L^{tn-\sum_{j\in
J}m_j-k}\cdot (\L-1)^{|J|}$$ hence $[B_t]\in FL\mathcal{N}$ and
$[B'_t]=[B_t]\L^k$. Therefore $h_{\infty}(C_k)=\pi_t^{-1}(B_t)$ is
$FL$-measurable and $\mu(C_k)=\mu(h_{\infty}(C_k))\cdot \L^k$.

Since $$\mathcal{L}(X)=\coprod_{k\in \Z_{\geq 0}\cup \infty}
h_{\infty}(C_k),$$ we have $[X]=\sum_{k\in \Z_{\geq
0}}\mu(h_{\infty}(C_k))=\sum_{k\in \Z_{\geq
0}}\mu(C_k)\L^{-k}=\int_{\mathcal{L}(Y)}\L^{-ord_t D}d\mu$.
\end{proof*}

\begin{definition}
Two smooth projective varieties $X_1, X_2\in VarFL$ are said to be
$FLK$-equivalent if there exists a smooth projective variety $Y\in
VarFL$ and two birational morphisms $\rho_1:Y \rightarrow X_1,
\rho_2:Y \rightarrow X_2$ such that
$\rho_1^*K_{X_1}=\rho_2^*K_{X_2}$ and the effective divisor
$D=K_Y-\rho^*_1K_{X_1}=\sum^r_{i=1}a_iD_i$ has simple normal
crossings and $D^0_J\in VarFL$ for any $J\subset \{1,..., r\}$.
\end{definition}

Directly from the change of variables formula, we get the following
result.
\begin{corollary}
If two smooth projective varieties $X_1, X_2\in VarFL$ are
$FLK$-equivalent, then $[X_1]=[X_2]$ in $FL\hat{\mathcal{N}}$.
\end{corollary}

We recall that each Lawson homology group $L_rH_n(X; \Q)$ has an
inductive limit of mixed Hodge structure (see \cite{FM}). And by a
result of Walker (see \cite{Walker}), the localization sequence of
Lawson homology groups is a sequence of inductive limit of mixed
Hodge structures. Since we are considering finite Lawson homology
groups, an inductive limit of mixed Hodge structure is just a mixed
Hodge structure.

Fix $n\in \Z_{\geq 0}$. For $X \in VarFL$, let $h_{p, q}(L_nH_i(X;
\C))$ be the dimension of the $(p, q)$-type Hodge component in
$L_nH_i(X; \C)$. We define the Lawson-Deligne-Hodge polynomial of
$X\in VarFL$ to be
$$F_nE(X):=\sum_{p, q} F_nE^{p, q}(X)u^pv^q$$
where
$$F_nE^{p, q}(X):=\sum_{i\geq 2n}(-1)^ih_{p, q}(L_nH_i(X; \C))$$
And we define the $F_rL$-Euler characteristic of $X$ to be
$$F_rL(X)=\sum_k(-1)^k dim L_rH_k(X; \Q)$$
From the localization sequence of Lawson homology, it is not
difficult to see that $F_rL, F_nE$ induces a group homomorphism from
$FL\overline{\mathcal{N}}$, the image of $FL\mathcal{N}$ in
$FL\hat{\mathcal{N}}$, to $\Z$ and $\Z[u, v]$ respectively.

\begin{corollary}
Two $FLK$-equivalent smooth projective varieties have the same
Lawson-Deligne-Hodge polynomial and the $F_rL$-Euler characteristic
for any $r\geq 0$.
\end{corollary}

\subsection{Higher Chow groups}
In the previous section, we use only the properties of localization
sequences of Lawson homology, the homotopy property and the
projective bundle theorem. Since there are analogous theorems for
higher Chow groups, we can play the same game for higher Chow
groups. Let $VarCH$ be the collection of all quasi-projective $X$
whose higher Chow groups $CH^r(X, n)$ are all finitely generated for
any $r, n$. Then from the localization sequence of higher Chow
groups:
$$\cdots \rightarrow CH^{q-d}(Z, p) \rightarrow CH^q(X, p)
\rightarrow CH^q(U, p) \rightarrow CH^{q-d}(Z, p-1) \rightarrow
\cdots$$ where $Z\subset X$ is a closed subvariety of codimension
$d$ and $U$ is its complement. Hence if $X, Z$ are in $VarCH$, then
$U$ is in $VarCH$. Then we can form the Grothendieck group
$K_0(VarCH)$ of $VarCH$. Similar to what we have done for Lawson
homology, we have some analogous results. We form $CH\mathcal{N}$,
$CH\hat{\mathcal{N}}$ and $CH\overline{\mathcal{N}}$ as their
analogs in Lawson homology.

\begin{definition}
Two smooth projective varieties $X_1, X_2\in VarCH$ are said to be
$CHK$-equivalent if there exists a smooth projective variety $Y\in
VarCH$ and two birational morphisms $\rho_1:Y \rightarrow X_1,
\rho_2:Y \rightarrow X_2$ such that
$\rho_1^*K_{X_1}=\rho_2^*K_{X_2}$ and the effective canonical
relative divisor $D=K_Y-\rho^*_1K_{X_1}=\sum^r_{i=1}a_iD_i$ has
simple normal crossings and $D^0_J\in VarCH$ for any $J\subset
\{1,..., r\}$.
\end{definition}

\begin{theorem}
Suppose that $X_1, X_2\in VarCH$ are two smooth projective varieties
which are $CHK$-equivalent. Then $[X_1]=[X_2]$ in
$CH\overline{\mathcal{N}}$.
\end{theorem}

\section{Projective bundle theorem} In \cite{F4}, Friedlander
proved a projective bundle theorem in morphic cohomology for smooth
normal quasi-projective varieties. In this section we prove a
projective bundle theorem of trivial bundles for all normal
quasi-projective varieties without assuming smoothness. Since we do
not have a Mayer-Vietoris sequence in morphic cohomology, the proof
is much more complicated than its counterpart in Lawson homology. We
are not sure if our approach may work for general bundles.

\begin{definition}
Suppose that $X$ is a normal quasi-projective variety and $W, Y$ are
projective varieties. Let $Z_k(W)(Y)$ be the subgroup of
$Z_{k+m}(W\times Y)$ consisting of algebraic cycles equidimensional
over $Y$ where $m$ is the dimension of $Y$ and $k$ is the dimension
of a fibre. The Chow variety $\mathscr{C}_{r, d}(W\times Y)$ of
$r$-dimensional algebraic cycles of degree $d$ of $W\times Y$ is a
projective variety and $\mathscr{M}or(X, \mathscr{C}_{r, d}(W\times
Y))$, the collection of all algebraic morphisms from $X$ to
$\mathscr{C}_{r, d}(W\times Y)$, is enrolled with the topology of
convergence with bounded degree (see \cite{FL3}). We define
$\mathscr{M}or(X, \mathscr{C}_{r, d}(W)(Y))$ to be the subspace of
$\mathscr{M}or(X, \mathscr{C}_{r, d}(W\times Y))$ consisting of
morphisms $f$ such that $f(x)\in \mathscr{C}_{r, d}(W)(Y)$ for all
$x\in X$ where $\mathscr{C}_{r, d}(W)(Y)$ is the collection of all
algebraic $r$-cycles of degree $d$ of $Y\times W$ equidimensional
over $Y$. Let
$$\mathscr{M}or(X,\mathscr{C}_r(W)(Y)):=\coprod_{d\geq
0}\mathscr{M}or(X, \mathscr{C}_{r, d}(W)(Y))$$ which is a
topological monoid and let
$$\mathscr{M}or(X, Z_r(W)(Y)):=[\mathscr{M}or(X,
\mathscr{C}_r(W)(Y))]^+$$ be its naive group completion. The
inclusion map $i:\mathscr{M}or(X, \mathscr{C}_r(W)(Y))
\hookrightarrow \mathscr{M}or(X, \mathscr{C}_{r+m}(W\times Y))$
induces a continuous homomorphism
$$\mathscr{D}':\mathscr{M}or(X, Z_k(W)(Y)) \rightarrow
\mathscr{M}or(X, Z_{k+m}(W\times Y))$$
\end{definition}

\begin{lemma}
Let $Y$ be a projective variety. From the Friedlander-Lawson moving
lemma (see \cite{FL2}), we have a map $\widetilde{\Psi}_t=(\Psi_{1,
t}, \Psi_{2, t}):\mathscr{C}_r(Y)\rightarrow \mathscr{C}_r(Y)\times
\mathscr{C}_r(Y)$ for $t\in I$ where $I$ is the unit interval $[0,
1]$ and $\mathscr{C}_r(Y)=\coprod_{d\geq 0}\mathscr{C}_{r, d}(Y)$ is
the Chow monoid of $r$-cycles of $Y$. The restriction of the maps
$$\Psi_{1, t}, \Psi_{2, t}:\mathscr{C}_{r, d}(Y) \rightarrow
\mathscr{C}_{r, d'}(Y)$$ are algebraic morphisms.
\end{lemma}

\begin{proof*}
The addition of effective cycles $C_1+C_2$ corresponding to the
product of their Chow forms $F_{C_1}F_{C_2}$. Hence the addition
$$+:\mathscr{C}_{r, d_1}(Y)\times \mathscr{C}_{r, d_2}(Y)
\rightarrow \mathscr{C}_{r, d_1+d_2}(Y)$$ is algebraic.

In \cite[Theorem 3.1]{FL2}, for an effective cycle $Z$ on $Y$,
$\Psi_{1, t}(Z), \Psi_{2, t}(Z)$ are the positive and negative parts
of $\Psi_{\underline{\textbf{N}}}-\Psi_{\underline{\textbf{N'}}}$
where
$$\Psi_{\underline{\textbf{N}}}=(-1)^{m+1}(M+1)R_{F^*}(Z)+\sum^m_{i=1}(-1)^i\pi^*_{F^i}\{\Theta_{\underline{\textbf{N}},
p}\{p_{F^i_*}(R_{F^{i-1}}\circ \cdots \circ R_{F^0}(Z))\}\}\bullet
Y.$$ Since all the maps involved are algebraic, the restriction of
$\Psi_{1, t}, \Psi_{2, t}$ to Chow varieties are algebraic.
\end{proof*}

\begin{theorem}
Suppose that $X$ is a normal quasi-projective variety and $W, Y$ are
smooth projective varieties. Then the map
$$\mathscr{D}':\mathscr{M}or(X, Z_r(W)(Y)) \rightarrow
\mathscr{M}or(X, Z_{r+m}(W\times Y))$$ is a homotopy equivalence.
\end{theorem}

\begin{proof*}
Let $\widetilde{\Psi}_t=(\Psi_{1, t}, \Psi_{2, t})$ be the map from
Friedlander-Lawson moving lemma where $t\in I $. By abuse of
notation, we define a map $\Psi_t:\mathscr{M}or(X, Z_r(W)(Y))
\rightarrow \mathscr{M}or(X, Z_r(W)(Y))$ by
$$\Psi_t(f):=\Psi_{1, t}(f)-\Psi_{2, t}(f)$$ and
$\Psi_{i, t}(f)(x):=\Psi_{i, t}(f(x))$ for $x\in X, i=1, 2$.

Let $$K_d:=\coprod_{d_1+d_2\leq d}\mathscr{M}or(X, \mathscr{C}_{r,
d_1}(W\times Y))\times \mathscr{M}or(X, \mathscr{C}_{r, d_2}(W\times
Y))/\sim,$$

$$K'_d:=\coprod_{d_1+d_2\leq d}\mathscr{M}or(X, \mathscr{C}_{r,
d_1}(W)(Y))\times \mathscr{M}or(X, \mathscr{C}_{r,
d_2}(W)(Y))/\sim$$ where $(f_1, g_1)\sim (f_2, g_2)$ if and only if
$f_1+g_2=f_2+g_1$.

The topology of $\mathscr{M}or(X, Z_{r+m}(W\times Y))$ is same as
the weak topology defined by the filtration
$$K_0\subset K_1\subset
K_2 \subset \cdots =\mathscr{M}or(X, Z_{r+m}(W\times Y))$$ and the
topology of $\mathscr{M}or(X, Z_r(W)(Y))$ is same as the weak
topology defined by the filtration
$$K'_0\subset K'_1\subset
K'_2 \subset \cdots =\mathscr{M}or(X, Z_{r+m}(W)(Y))$$ From
\cite[Lemma 2.3]{T1}, these two filtrations are locally compact. Let
$\phi_{e, t}, \phi'_{e, t}$ be the restriction of $\Psi_{1, t},
\Psi_{2, t}$ to $K_e$ and $K'_e$ respectively. Let
$\lambda_e:K_e\times \{1\} \rightarrow \mathscr{M}or(X, Z_r(W)(X))$
be $\phi_{e, 1}$. Then we have the following diagrams:
$$\xymatrix{K'_e\times I \ar[r]^-{\phi'_e} \ar[d]_{\mathscr{D}'\times Id} & \mathscr{M}or(X,
Z_r(W)(Y)) \ar[d]^{\mathscr{D}'} \\
K_e\times I \ar[r]^-{\phi_e} & \mathscr{M}or(X, Z_{r+m}(W\times
Y))}$$

$$\xymatrix{K'_e\times \{1\} \ar[d]_{\mathscr{D}'} \hookrightarrow K'_e\times I \ar[r]^-{\phi'_e}
& \mathscr{M}or(X, Z_r(W)(Y)) \ar[d]^{\mathscr{D}'}\\
K_e\times \{1\} \hookrightarrow K_e\times I \ar[ru]^{\lambda_e}
\ar[r]^-{\phi_e} & \mathscr{M}or(X, Z_{r+m}(W\times Y))}
$$

Then by \cite[Lemma 5.2]{FL3}, $\mathscr{D}'$ is a weak homotopy
equivalence. Since $\mathscr{M}or(X, Z_r(W)(Y))$ and
$\mathscr{M}or(X, Z_{r+m}(W\times Y))$ have the homotopy type of a
CW-complex, by Whitehead theorem, $\mathscr{D}'$ is a homotopy
equivalence.
\end{proof*}

\begin{proposition}
Suppose that $X$ is a normal quasi-projective variety and $W, Y$ are
projective varieties. Then $\mathscr{M}or(X\times Y, Z_r(W))$ is
isomorphic as a topological group to $\mathscr{M}or(X, Z_r(W)(Y))$.
\end{proposition}

\begin{proof*}
There is a natural bijection
$$\psi:\mathscr{M}or(X\times Y, \mathscr{C}_r(W)) \rightarrow \mathscr{M}or(X,
\mathscr{C}_r(W)(Y))$$ defined by $\psi(f)(x)(y):=f(x, y)$. These
two spaces are obviously homeomorphic under the topology of
convergence with bounded degree and $\psi$ is monoid isomorphism.
Hence we complete the proof.
\end{proof*}

Consider the topology of convergence with bounded degree, we have
the following fact.
\begin{lemma}
For a normal quasi-projective variety $X$ and a projective variety
$Y$, let $\mathscr{M}or(X, Z_{r_1}(Y)\times \cdots \times
Z_{r_k}(Y))$ be the topological naive group completion of
$\mathscr{M}or(X, \mathscr{C}_{r_1}(Y)\times
\mathscr{C}_{r_2}(Y)\times \cdots \times \mathscr{C}_{r_k}(Y))$.
Then there is an isomorphism of topological groups:
$$\mathscr{M}or(X, Z_{r_1}(Y)\times Z_{r_2}(Y)\times \cdots \times
Z_{r_k}(Y))\cong\mathscr{M}or(X, Z_{r_1}(Y))\times \mathscr{M}or(X,
Z_{r_2}(Y))\times \cdots \times \mathscr{M}or(X, Z_{r_k}(Y))
$$
\end{lemma}

\begin{definition}
For normal quasi-projective varieties $X, U$, if $U=Y-Z$ where $Y,
Z$ are projective varieties and $Z$ is a subvariety of $Y$, then we
define
$$\mathscr{M}or(X, Z_i(U)):=\frac{\mathscr{M}or(X, Z_i(Y))}{\mathscr{M}or(X,
Z_i(Z))}$$
\end{definition}

\begin{proposition}(localization sequence)
For $X, U, Y, Z$ as above. There is a localization sequence:
$$\cdots \rightarrow \pi_k\mathscr{M}or(X, Z_i(Z)) \rightarrow
\pi_k\mathscr{M}or(X, Z_i(Y)) \rightarrow \pi_k\mathscr{M}or(X,
Z_i(U)) \rightarrow \pi_{k-1}\mathscr{M}or(X, Z_i(Z)) \rightarrow
\cdots $$
\end{proposition}

\begin{definition}
For a normal quasi-projective variety $X$ and a projective variety
$Y$, we define
$$\mathscr{M}or(X, Z^t(Y)):=\frac{\mathscr{M}or(X, Z_0(\P^t)(Y))}{\mathscr{M}or(X, Z_0(\P^{t-1})(Y))}$$
with the quotient topology.
\end{definition}

\begin{theorem}
For $X$ a normal quasi-projective variety and $Y$ a smooth
projective variety,
\begin{enumerate}
\item $\mathscr{M}or(X, Z_r(Y\times \A^t))$ is homotopy equivalent
to $\mathscr{M}or(X, Z_{r-t}(Y))$ if $r\geq t$.

\item (Duality)
$\mathscr{M}or(X, Z^t(Y))$ is homotopy equivalent to
$\mathscr{M}or(X, Z_{n-t}(Y))$.
\end{enumerate}
\end{theorem}

\begin{proof*}
\begin{enumerate}
\item
By \cite[Proposition 3.7]{F4}, the suspension $\susp_*:
\mathscr{M}or(X, Z_r(Y)) \rightarrow \mathscr{M}or(X, Z_{r+1}(\susp
Y))$ is a homotopy equivalence. Let $L$ be the restriction of the
hyperplane line bundle $\mathscr{O}(1)$ of $\P^N$ to $Y$ and observe
that $L=\susp Y-\{\infty\}$. Hence $\mathscr{M}or(X, Z_r(Y))$ is
homotopy equivalent to $\mathscr{M}or(X, Z_{r+1}(L))$. Now the
result follows from an induction on dimension of $Y$ and the
localization sequence.

\item
Consider the following commutative diagram:
$$\xymatrix{\mathscr{M}or(X, Z_0(\P^{t-1})(Y)) \ar[r] \ar[d]^{\mathscr{D}}&
\mathscr{M}or(X, Z_0(\P^t)(Y)) \ar[r] \ar[d]^{\mathscr{D}} & \mathscr{M}or(X, Z^t(Y)) \ar[d]\\
\mathscr{M}or(X, Z_n(Y\times \P^{n-1})) \ar[r] & \mathscr{M}or(X,
Z_n(Y\times \P^t)) \ar[r] & \mathscr{M}or(X, Z_n(Y\times \A^t))}$$
where $Z^t(Y)=\frac{Z_0(\P^t)(Y)}{Z_0(\P^{t-1})(Y)}$. From the
induced long exact sequence of homotopy groups and the result above,
we see that $\mathscr{M}or(X, Z^t(Y))$ is homotopy equivalent to
$\mathscr{M}or(X, Z_n(Y\times \A^t))\cong \mathscr{M}or(X,
Z_{n-t}(Y))$.
\end{enumerate}
\end{proof*}

\begin{theorem}
For a projective variety $Y$, there is a splitting
$$\xi:\mathscr{M}or(X, Z_0(\P^t)(Y))\rightarrow \mathscr{M}or(X, Z^t(Y))\times
\mathscr{M}or(X, Z^{t-1}(Y))\times \cdots \times \mathscr{M}or(X,
Z^0(Y))$$ which is a homotopy equivalence.
\end{theorem}

\begin{proof*}
By the construction in the proof of the splitting theorem of Lawson
and Friedlander (\cite[Theorem 2.10]{FL1}), there is a projection
map
$$p:Z_0(\P^t)(Y) \rightarrow Z_0(\P^t)(Y)\times
Z_0(\P^{t-1})(Y)\times \cdots \times Z_0(\P^0)(Y)$$ Write $p=(p_t,
p_{t-1}, ..., p_0)$ and for $f\in \mathscr{M}or(X, Z_0(\P^t)(Y))$,
define
$$p_i(f)(x):=p_i(f(x))$$
Then we get a map
$$\xi^t:\mathscr{M}or(X, Z_0(\P^t)(Y))\rightarrow \mathscr{M}or(X,
Z^t(Y))\times \mathscr{M}or(X, Z^{t-1}(Y))\times \cdots \times
\mathscr{M}or(X, Z^0(Y))$$ defined by
$$\xi^t(f)=(p_t(f)+\mathscr{M}or(X, Z_0(\P^{t-1})(Y)), p_{t-1}(f)+\mathscr{M}or(X,
Z_0(\P^{t-2})(Y)), ..., p_0(f)+\mathscr{M}or(X, Z_0(\P^0)(Y))$$

We are going to show this map is a homotopy equivalence. We prove by
induction on $t$. When $t=0$, this follows from definition. Assume
that we have the splitting for $t-1$. Consider the following
commutative diagram
$$\xymatrix{\mathscr{M}or(X, Z_0(\P^{t-1})(Y)) \ar[d] \ar[rr]^-{\xi^{t-1}} && \mathscr{M}or(X, Z^{t-1}(Y))
\times \cdots \times \mathscr{M}or(X, Z^0(Y)) \ar[d] \\
\mathscr{M}or(X, Z_0(\P^t)(Y)) \ar[d] \ar[rr]^-{\xi^t}
&&\mathscr{M}or(X, Z^t(Y))\times \mathscr{M}or(X, Z^{t-1}(Y))\times
\cdots \times
\mathscr{M}or(X, Z^0(Y)) \ar[d]\\
\mathscr{M}or(X, Z^t(Y)) \ar[rr]^{=}&& \mathscr{M}or(X, Z^t(Y))}$$
Follow from the long exact sequences of homotopy groups induced from
the vertical rows, we see that $\xi^t$ is a weak homotopy
equivalence. Since all spaces have the homotopy type of a
CW-complex, by the Whitehead theorem, $\xi^t$ is a homotopy
equivalence.
\end{proof*}

Combining the splitting and the duality theorem, we get the
following splitting.

\begin{corollary}
For a normal quasi-projective variety $X$ and $Y$ a smooth
projective of dimension $n$, if $t\leq n$, there is a splitting
$$\xi': \mathscr{M}or(X, Z_0(\P^t)(Y))\cong \mathscr{M}or(X,
Z_{n-t}(Y))\times \mathscr{M}or(X, Z_{n-t+1}(Y))\times \cdots \times
\mathscr{M}or(X, Z_n(Y))$$ which is a homotopy equivalence.
\end{corollary}

\begin{theorem}
For a normal quasi-projective variety $X$, there is a homotopy
equivalence $\overline{\eta}:Z^t(X\times \P^e) \cong
\bigoplus^t_{i=0}Z^{t-i}(X)$ for $e\geq t$.
\end{theorem}

\begin{proof*}
By Lawson suspension theorem (see \cite[Proposition 3.7]{F4}), there
is a homotopy equivalent $\susp_*:\mathscr{M}or(X,
Z_{e-t+i}(\P^e))\cong \mathscr{M}or(X, Z_0(\P^{t-i}))$, We have a
homotopy equivalence
$$\eta: \mathscr{M}or(X\times \P^e,
Z_0(\P^t))\overset{\psi}{\longrightarrow} \mathscr{M}or(X,
Z_0(\P^t)(\P^e))\overset{\xi}{\longrightarrow}\bigoplus^t_{i=0}\mathscr{M}or(X,
Z^{t-i}(\P^e))\overset{\mathscr{D}}{\longrightarrow}\bigoplus^t_{i=0}\mathscr{M}or(X,
Z_0(\P^{t-i}))$$

Consider the following commutative diagram:
$$\xymatrix{\mathscr{M}or(X\times \P^e, Z_0(\P^{t-1})) \ar[d]^{\eta} \ar[r]
& \mathscr{M}or(X\times \P^e, Z_0(\P^t)) \ar[d]^{\eta} \ar[r] &
Z^t(X\times
\P^e)\ar[d]^{\overline{\eta}}\\
\bigoplus^{t-1}_{i=0}\mathscr{M}or(X, Z_0(\P^{t-1-i})) \ar[r] &
\bigoplus^t_{i=0}\mathscr{M}or(X, Z_0(\P^{t-i})) \ar[r] &
\bigoplus^t_{i=0}Z^{t-i}(X)}$$ From the long exact sequences of
homotopy groups induced by the horizontal rows, we see that
$\overline{\eta}$ is a weak homotopy equivalence. But all these
spaces have the homotopy type of CW-complexes, hence
$\overline{\eta}$ is a homotopy equivalence.
\end{proof*}

The morphic cohomology groups are known only for very few cases of
smooth varieties, and almost nothing about singular varieties. As an
application of above result, we calculate the morphic cohomology
groups of two singular surfaces.

\begin{example}
One of the main tools we use is \cite[Theorem 9.1]{FL1} which says
that there is a fibration $Z^1(X) \rightarrow Pic(X)$ with homotopy
fibre $K(\Z, 2)$ for projective variety $X$.
\begin{enumerate}
\item
Let $C_1:y^2Z=x^2(x+z)$ be the curve in $\P^2$ which has a node at
$[0:0:1]$. The Picard group of $C_1$ is $Pic(C_1)\cong \C^*\times
\Z$. Therefore we know that $L^1H^k(C_1\times \P^1)$ is 0 for $k>
2$. We list other cases in the following table.

\begin{tabular}{|c|c|c|c|c|}
  \hline
  $k$& $\pi_kPic(C_1)$ & $\pi_kZ^1(C_1)$ & $\pi_kZ^0(C_1)$ & $L^1H^k(C_1\times \P^1)$\\
\hline
  0 & $\Z$ & $\Z$ & $\Z$ & $\Z$  \\
\hline
  1 & $\Z$ & $\Z$  & 0 & $\Z$ \\
\hline
  2 &   0 &  $\Z$  & 0 & $\Z\oplus \Z$ \\
  \hline
\end{tabular}

\item
Let $C_2:y^2z=x^3$ be the curve in $\P^2$ which has a cusp at
$(0:0:1)$. The Picard group of $C_2$ is $Pic(C_2)\cong \C\times \Z$.
Hence $L^1H^k(C_2\times \P^1)$ is 0 for $k> 2$. We list other cases
in the following table.

\begin{tabular}{|c|c|c|c|c|c|c|c|}
  \hline
  $k$& $\pi_kPic(C_2)$ &$\pi_kZ^1(C_2)$ & $\pi_kZ^0(C_2)$ & $L^1H^k(C_2\times \P^1)$ \\
\hline
  0 & $\Z$ & $\Z$ & $\Z$ & $\Z$  \\
\hline
  1 & 0 & 0  & 0 & 0 \\
\hline
  2 &   0 &  $\Z$  & 0 & $\Z\oplus \Z$ \\
  \hline
\end{tabular}

\end{enumerate}
\end{example}

\begin{acknowledgements}
The author thanks Chin-Lung Wang for a question he asked, Willem
Veys and Francois Loeser for answering some of his trivial questions
and the referee for his very nice comments and suggestions. He
thanks the National Center for Theoretical Sciences in Hsinchu,
Taiwan for financial support.
\end{acknowledgements}

\end{document}